\documentclass{article}
\usepackage{graphicx} 
\usepackage[top=1in,bottom=1in,left=1in,right=1in]{geometry}
\usepackage{amssymb}
\usepackage{amsmath}
\usepackage{amsbsy}
\usepackage{mathtools}
\usepackage{xcolor}
\usepackage{hyperref}
\usepackage{graphicx}
\usepackage{xcolor}
\usepackage{caption}
\usepackage{subcaption}
\usepackage{float}
\usepackage{multirow}
\usepackage{adjustbox}
\usepackage{relsize}
\usepackage{enumitem}

\newcommand{\update}[1]{\textcolor{black}{#1}} 

\DeclareMathOperator*{\argmin}{arg\,min}
\def\sol{S}
\def\initCond{J}
\def\R{\mathbb{R}}
\def\Rn{\R^n}
\def\LOvar{u}
\def\Hopfvar{p}

\title{Recent Advances in Numerical Solutions for Hamilton-Jacobi PDEs}
\author{Tingwei Meng\textsuperscript{1} \and 
        Siting Liu\textsuperscript{2} \and 
        Samy Wu Fung\textsuperscript{3} \and 
        Stanley Osher\textsuperscript{1}}
\date{}
\begin{document}
\maketitle

\footnotetext[1]{Department of Mathematics, University of California, Los Angeles}
\footnotetext[2]{Department of Mathematics, University of California, Riverside}
\footnotetext[3]{Department of Applied Mathematics and Statistics, Colorado School of Mines}

\begin{abstract}
    Hamilton-Jacobi partial differential equations (HJ PDEs) play a central role in many applications such as economics, physics, and engineering. 
    These equations describe the evolution of a value function which encodes valuable information about the system, such as action, cost, or level sets of a dynamic process. Their importance lies in their ability to model diverse phenomena, ranging from the propagation of fronts in computational physics to optimal decision-making in control systems.
    This paper provides a review of some recent advances in numerical methods to address challenges such as high-dimensionality, nonlinearity, and computational efficiency. By examining these developments, this paper sheds light on important techniques and emerging directions in the numerical solution of HJ PDEs.
\end{abstract}

\section{Introduction}
Hamilton-Jacobi partial differential equations (HJ PDEs) are foundational mathematical tools that arise in diverse areas such as optimal control, geometric optics, and computational physics~\update{\cite{evans2022partial, bardi1997optimal, crandall1983viscosity,crandall1984some}}. These equations are commonly written in the form  
\begin{equation}\label{eqt:HJ_general}
\begin{dcases}
\partial_t {\sol}(x,t) + H(\nabla_x {\sol}(x,t), x, t) = 0, & x \in \Omega, t > 0, \\
\sol(x,0) = \initCond(x), & x \in \Omega,
\end{dcases}
\end{equation}
where $\sol(x, t)$ is the value function, $H$ is the Hamiltonian encoding the system's dynamics, $\initCond(x)$ represents the initial condition, and $\Omega \subseteq \Rn$ is the spatial domain. If $\Omega$ is a bounded domain, certain boundary conditions are imposed on $S$.
If $\Omega$ is the whole space $\Rn$, the value $S(y,T-t)$ is the optimal value of the following optimal control problem: 
\begin{equation}\label{eqt:optimal_control}
\min\left\{\int_t^T L(u(s), x(s), s) ds + J(x(T)) \colon \dot x(s) = f(u(s),x(s),s) \; \forall s\in [t,T], x(t) = y\right\},
\end{equation}
where $u$ is the control function taking value in $\R^m$. In these two problems, the running cost $L$, the dynamical source term $f$, and the Hamiltonian $H$ are connected by
\begin{equation*}
H(p,x,t) = \sup_{u\in \R^m} \{-\langle p, f(u,x,T-t) \rangle - L(u, x, T-t) \}.
\end{equation*}

HJ PDEs are notable for their ability to capture complex dynamics in high-dimensional and nonlinear systems~\cite{fleming2006controlled, peng1999pde}. They are central to modeling phenomena such as the propagation of wavefronts~\cite{osher1988fronts}, {the computation of geodesics~\cite{bazanski1989hamilton, witzany2019hamilton}}, and decision-making in optimal control problems~\cite{bardi1997optimal, evans2022partial}. 
Despite their theoretical significance and wide applicability, solving HJ PDEs efficiently, especially in high dimensions, remains a challenging task, \update{as} traditional grid-based methods face the curse of dimensionality~\cite{beck2023overview, onken2022neural}.
{
In addition to their role in control theory, HJ PDEs capture fundamental principles in physics and geometry such as the Eikonal equation, which describes the travel-time of a wavefront in a medium with refractive index. In classical mechanics, $S(x,t)$ can be viewed as encoding the action integral, while the Hamiltonian $H$ represents the system’s total energy. 
Through this lens, the solution of the HJ PDE governs not only optimal trajectories in control but also the evolution of rays in optics and the variational structure of mechanics.
}

Recent years have seen significant progress in the development of numerical methods for \update{approximating viscosity solutions to HJ PDEs as defined by~\cite{crandall1983viscosity}}. These advances have been driven by the need for computationally efficient, accurate, and scalable algorithms to tackle high-dimensional problems and non-smooth solutions~\cite{darbon2016algorithms}. Techniques leveraging variational principles~\cite{chow2019algorithm}, level-set methods~\cite{osher1988fronts}, and modern machine learning~\cite{lin2021alternating, onken2021neural, ruthotto2020machine} have opened new avenues for addressing these challenges.  

This paper reviews some of the recent advances in numerical methods for HJ PDEs, highlighting their theoretical underpinnings, practical implementations, and limitations. 
We review a series of approaches including grid-based methods (Section~\ref{sec: grid_based}), methods based on representation formulas (Section~\ref{sec: representation_formulas}), Monte-Carlo-based methods using Laplace's Method (Section~\ref{sec: Laplace_approximations}), and Deep Learning methods (Section~\ref{sec: deep_learning}). While not comprehensive, our review of these developments aims to provide insights into the current numerical landscape and to identify promising directions for future research.


\section{Grid-based Methods}
\label{sec: grid_based}
Grid-based numerical methods have been the cornerstone of computational approaches for solving HJ PDEs~\cite{crandal1984two, osher1988fronts, osher1991high, liu1994weighted, abgrall1996numerical, shu2007high}.  Recent advances have focused on improving accuracy, computational efficiency, and handling of discontinuities in the solution gradients. \update{In addition to finite-difference/finite-volume discretizations,  semi-Lagrangian schemes provide an alternative grid-based route by discretizing the characteristic control problem and are especially effective for optimal-control HJ equations~\cite{falcone2013semi}.}  This section provides an overview of some key developments in grid-based methods over the past decade. For a more comprehensive review of grid-based numerical methods for HJ PDEs, we refer to~\cite{falcone2016numerical}.

\subsection{High-Order ENO and WENO Schemes}
While monotone schemes~\cite{crandal1984two, osher1988fronts, osher1991high} are inherently limited to first-order accuracy, high-order Essentially Non-Oscillatory (ENO) schemes~\cite{shu1988efficient} and their Weighted ENO (WENO) variants~\cite{liu1994weighted} have been developed to achieve uniform high-order accuracy in smooth regions while robustly capturing discontinuities in derivatives.
In the context of one-dimensional finite difference methods for solving HJ PDEs~\eqref{eqt:HJ_general}, for simplicity consider the case where $H(\nabla_x {\sol}(x,t), x, t) = H(\nabla_x {\sol}(x,t))$, the ENO or WENO method takes the form
\update{
\begin{align*}
\dfrac{d}{dt} S_i + \hat{H}\left(S_x^-(x_i,t),S_x^+(x_i,t)\right) = 0,   
\end{align*}}
where $S_i$ approximates the value of $S(x_i,t)$ \update{ and $\hat{H}(\cdot, \cdot)$ is a monotone numerical Hamiltonian, i.e., a consistent approximation of $ \hat{H}$ that is nondecreasing in its first argument and nonincreasing in its second.} \update{Here $S_x^-(x_i,t)$ and $S_x^+(x_i,t)$ denote, respectively, the left- and right-biased approximations of the derivative $S_x(x_i,t)$ obtained from ENO or WENO reconstructions.} \update{For the ENO scheme, the approximation of $S_x^-(x_i,t)$ starts with a stencil containing only ${S_{i-1}, S_i}$; the stencil is then adaptively extended by adding points that yield the smoothest local reconstruction—thereby avoiding interpolation across discontinuities—until the desired order is achieved, at which point a single optimally chosen stencil provides the approximation of $S_x^-(x_i,t)$. The computation of $S_x^+(x_i,t)$ proceeds similarly.} As for \update{the} WENO scheme, it forms a convex combination  of candidate approximations computed on several stencils. A typical approximation to the spatial derivatives can be written as:
\begin{equation*}
    S_x^-(x_i,t) = \sum_{k=0}^{r-1} \omega_k \hat{\sol}_{x_i}^k, \quad \omega_k = \frac{\alpha_k}{\sum_{l=0}^{r-1} \alpha_l}, \quad \alpha_k = \frac{d_k}{(\epsilon + \beta_k)^p},
\end{equation*}
where $\hat{\sol}_{x_i}^k$ is the $k$-th stencil approximation of $S_x^-(x_i,t)$, $\omega_k$ are non-linear weights satisfying $\sum_{k=0}^{r-1} \omega_k = 1$. The coefficients $d_k$ are the linear weights that form the approximation $\sum_{k=0}^{r-1} \omega_k \hat{\sol}_{x_i}^k$ as the $r$-th degree polynomial reconstruction over the full stencil. $\beta_k$ are smoothness indicators that measure the local variation of the solution, and $\epsilon$ is a small positive number to avoid division by zero. Based on the spatial discretizations discussed above, time \update{discretizations} can be performed using either explicit total variation diminishing Runge–Kutta methods or multi-step schemes.
Hermite polynomials are used to develop Hermite-type WENO schemes for solving HJ PDEs~\cite{qiu2005hermite,zheng2016directly}.

These methods have undergone continuous refinements to improve their performance, robustness, and computational efficiency.
Significant advances in WENO methodology include the development of exponential polynomial-based schemes by \cite{kim2021third}, achieving third-order accuracy with improved handling of discontinuities. In the same year, \cite{abedian2021rbf} introduced innovative RBF-ENO/WENO schemes combined with Lax-Wendroff type time discretizations, demonstrating superior performance for problems with complex solution structures.
Recent efforts have emphasized adaptivity and computational efficiency. In~\cite{abedian2023weno}, a fifth-order WENO scheme based on Legendre polynomials is constructed for simulating HJ equations in a finite difference framework. Most recently, \cite{zhang2025new} introduced a third-order, fifth-order, and seventh-order finite difference ghost multi-resolution WENO (GMR-WENO) schemes that only utilize the information defined on one four-point, one six-point, or one eight-point spatial stencil for designing high-order spatial approximations without introducing any other smaller stencils.

\subsection{Discontinuous Galerkin Methods}
\update{
The discontinuous Galerkin (DG) finite element method is a flexible, high-order approach that uses discontinuous piecewise polynomial spaces for both numerical solution and test functions~\cite{cockburn1991runge}. For Hamilton–Jacobi equations, two main DG strategies have been developed.  Early formulations evolve the spatial derivative of the solution and then recover the potential by local integration or reconstruction~\cite{hu1999discontinuous,li2005reinterpretation}.
Subsequent works apply DG directly to the solution variable, using one-sided derivative reconstructions at element interfaces together with the monotone numerical Hamiltonian (as defined above)~\cite{cheng2007discontinuous,li2010central,CHENG2014134}. Local DG (LDG) variants~\cite{yan2011local}, adaptive sparse-grid LDG for a high-dimensional setup~\cite{GUO2021110294}, the alternating-evolution DG method achieving high-order accuracy~\cite{liu2014alternating}, and arbitrary Lagrangian–Eulerian local DG on moving meshes~\cite{klingenberg2017arbitrary} further enhance accuracy, robustness, and efficiency.}

\section{Solutions Based on Representation Formulas}
\label{sec: representation_formulas}
For certain classes of HJ PDEs, viscosity solutions can be efficiently computed using representation formulas. A key advantage of this approach is that it transforms the problem of solving a PDE into an optimization problem, which helps mitigate the curse of dimensionality. In this section, we focus on HJ PDEs defined on the entire space \( \Omega = \mathbb{R}^n \) and explore how these representation formulas leverage intrinsic problem structures, such as constant or piecewise constant velocities or momenta, to facilitate computation.

\begin{figure}[htbp]
    \centering
    \begin{subfigure}{0.24\textwidth}
        \centering \includegraphics[width=\textwidth]{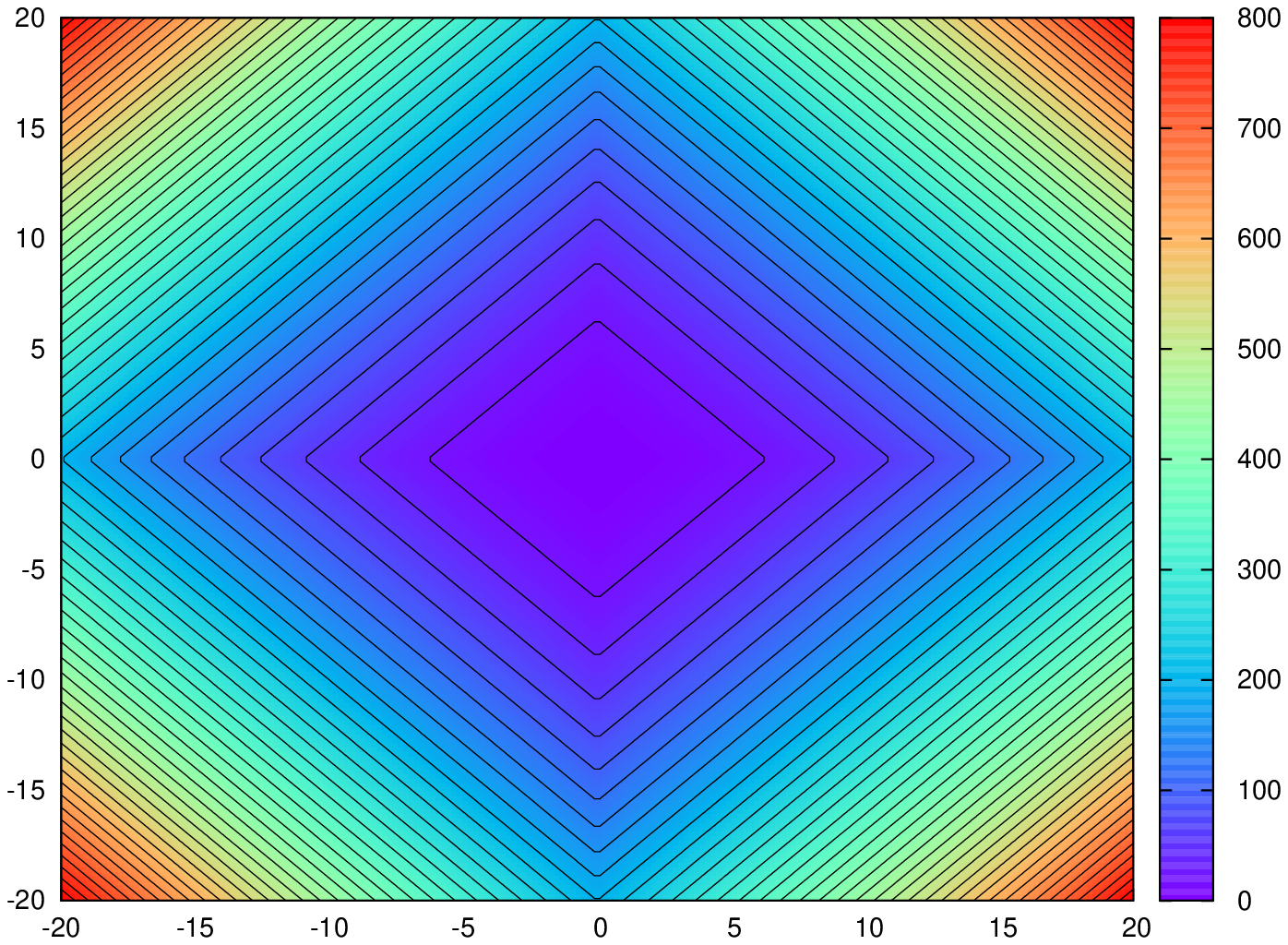}
    \end{subfigure}
    \begin{subfigure}{0.24\textwidth}
        \centering \includegraphics[width=\textwidth]{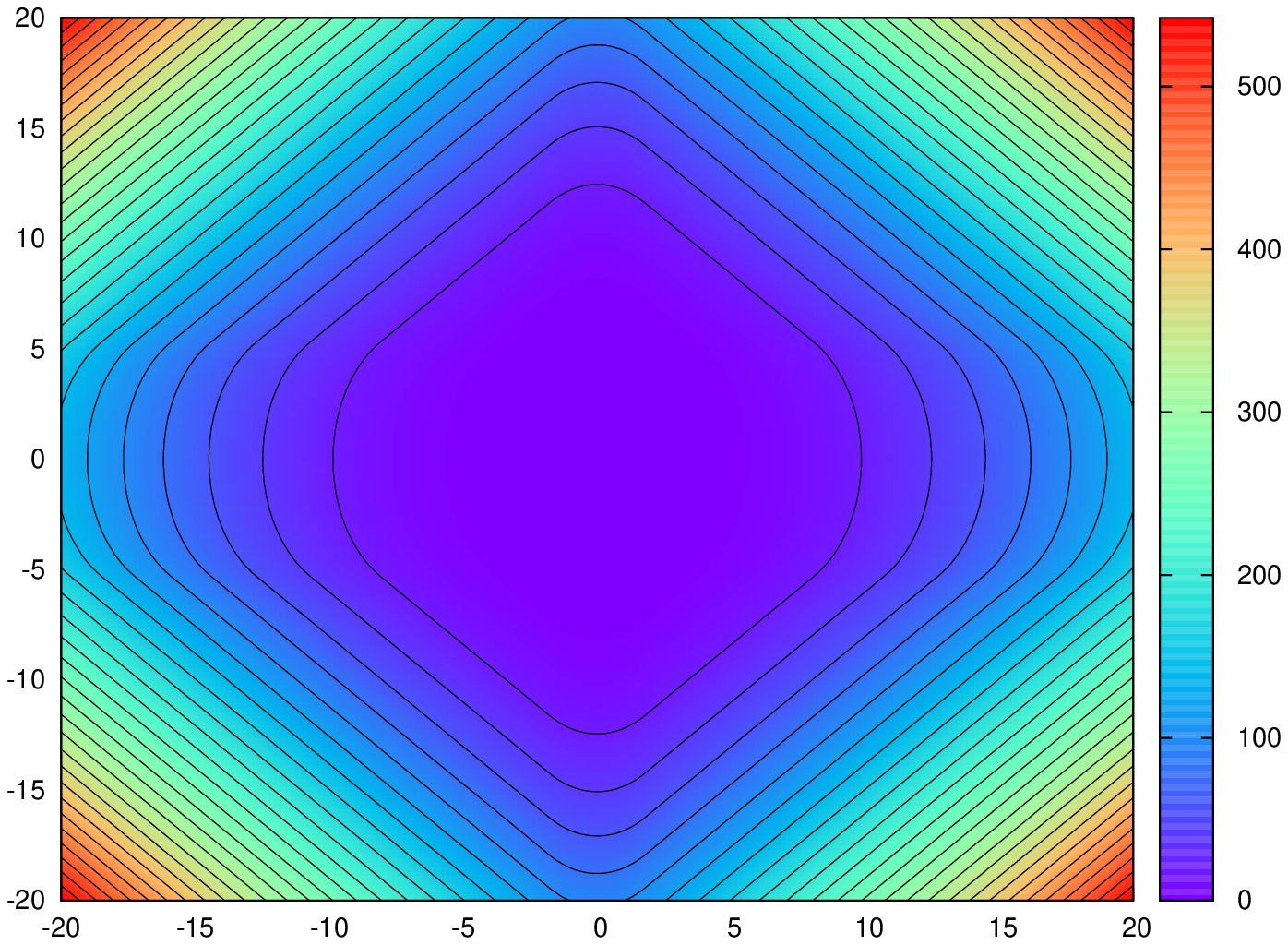}
    \end{subfigure}
    \begin{subfigure}{0.24\textwidth}
        \centering \includegraphics[width=\textwidth]{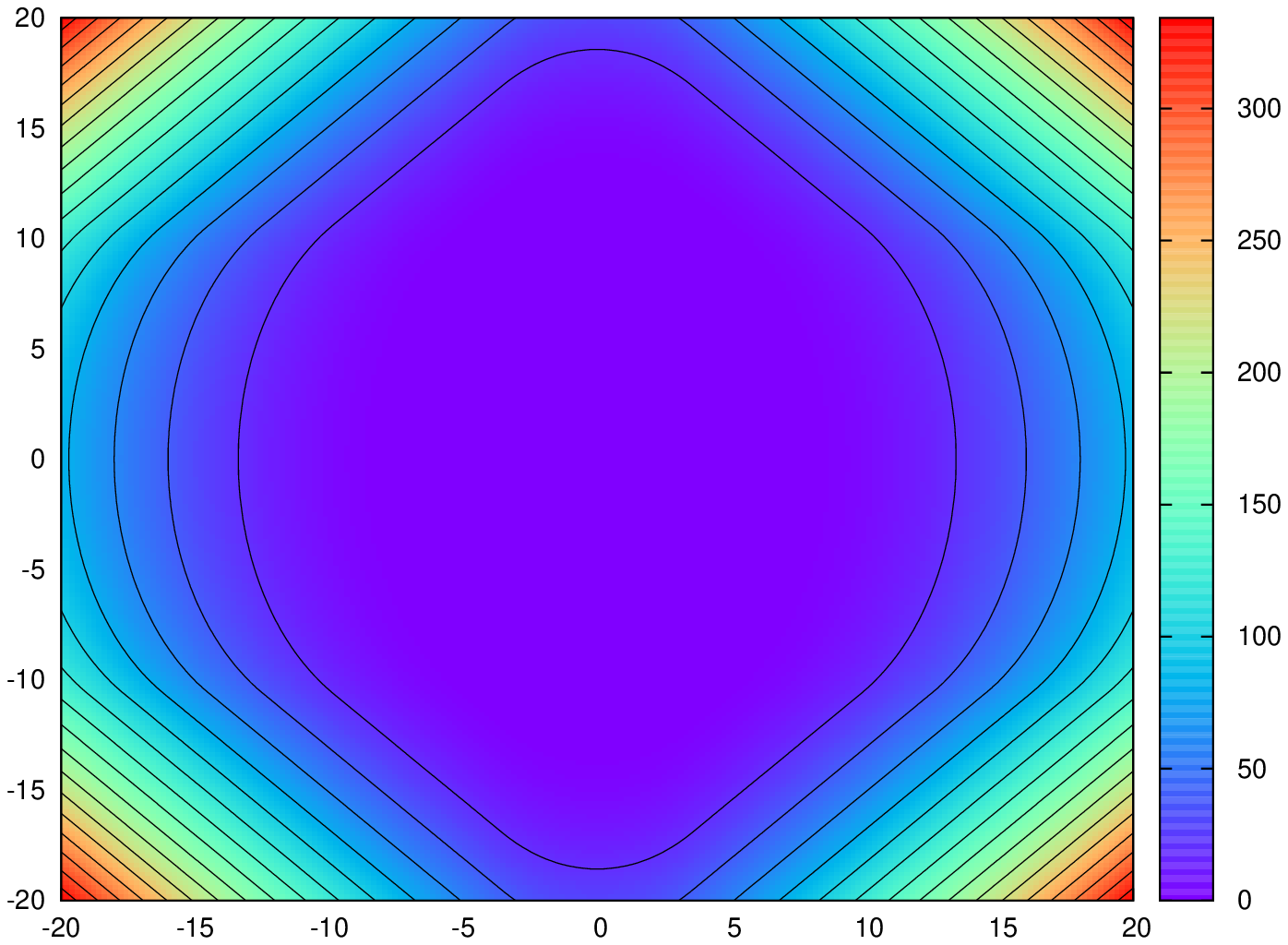}
    \end{subfigure}
    \begin{subfigure}{0.24\textwidth}
        \centering \includegraphics[width=\textwidth]{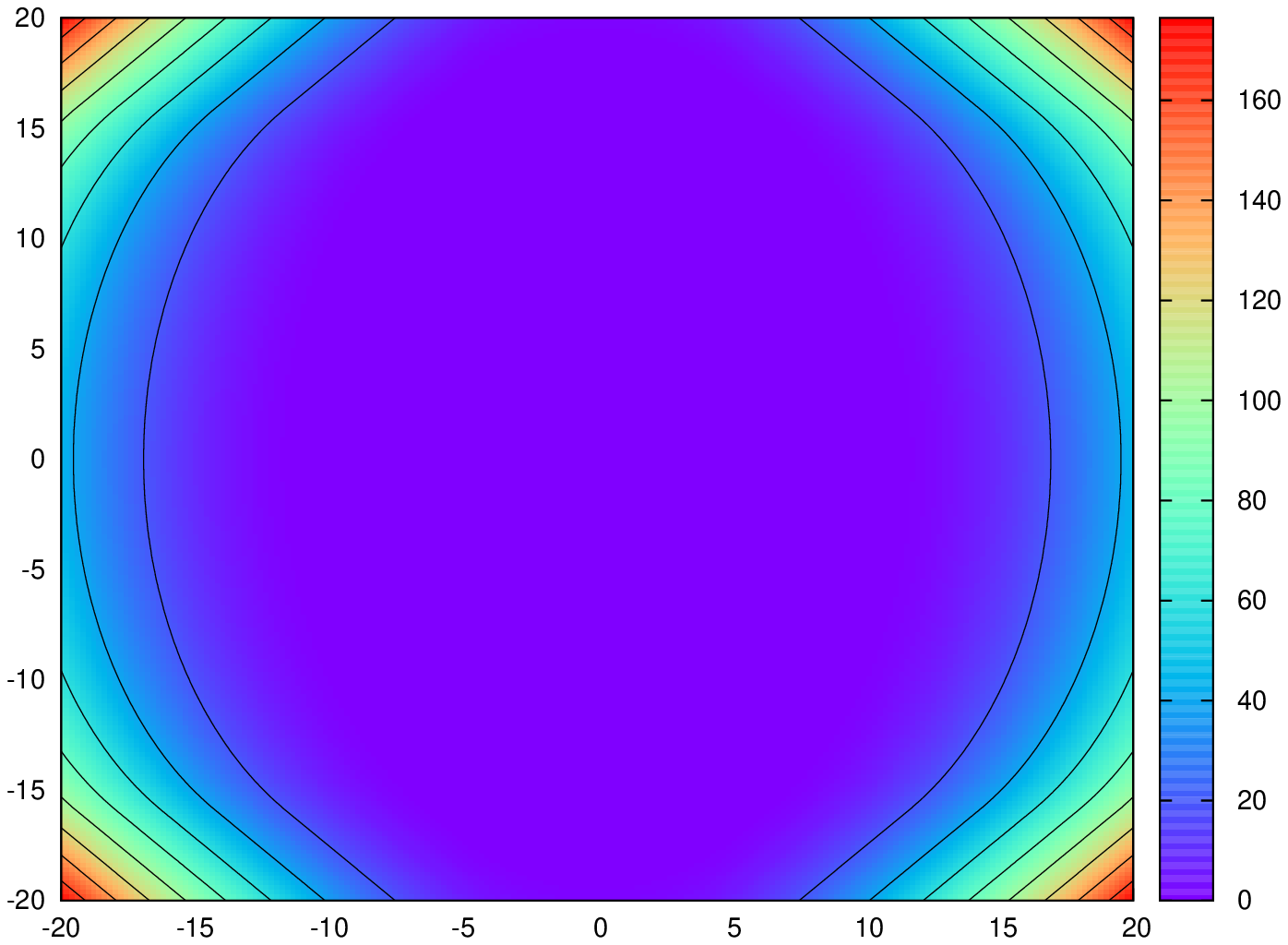}
    \end{subfigure}
    \caption{Evaluation of the solution $S((x_1,x_2,0,0,0,0,0,0)^T,t)$ of
    the HJ-PDE with initial data  $J=\frac{1}{2}\|\cdot\|_1^2$
    and Hamiltonian $H=\|\cdot\|_D$, for $(x_1,x_2) \,\in\,
    [-20,20]^2$ for   different times $t=0, 5,10,15$.
    }
    \label{fig:ch_line_results}
\end{figure}

Firstly, we consider the case when the Hamiltonian only depends on the momentum $\nabla_x S$. In this case, the viscosity solution has two representation formulas, namely Lax-Oleinik and Hopf formula~\cite{bardi1998hopf}. \update{We denote by $f^*$ the 
Fenchel–Legendre transform of a convex lower-semicontinuous function $f$, defined by $f^*(p) = \sup_{x\in\Rn} \{\langle x,p\rangle - f(x)\}$.} When the Hamiltonian $H$ is convex, the solution is given by the Lax-Oleinik formula
\begin{equation}
\sol(x,t) = \inf_{\LOvar\in \Rn} \left\{\initCond(\LOvar) + tH^*\left(\frac{x-\LOvar}{t}\right)\right\}.
\label{eq: lax_oleinik}
\end{equation}
On the other hand, when the initial condition $\initCond$ is convex, the solution is given by the Hopf formula
\begin{equation}\label{eqt:Hopf}
\sol(x,t) = \sup_{\Hopfvar\in \Rn} \left\{\langle x,\Hopfvar\rangle - tH(\Hopfvar) - \initCond^*(\Hopfvar)\right\}.
\end{equation}
These formulas can be used to compute the viscosity solution at a given point $(x,t)$ by solving an optimization problem. In this way, the algorithm can mitigate the curse of dimensionality. Moreover, the value on each point is independent from each other, making the computation easy to \update{parallelize}. When the Hamiltonian is 1-homogeneous,~\cite{darbon2016algorithms} designed fast solver by using the split Bregman iterative method to solve these two representation formulas. 
Later, in~\cite{lee2017revisiting}, the authors use these formulas to compute the redistancing problem, which, after using the level set method, can be solved by solving an Eikonal equation. {A numerical result is shown in Fig.~\ref{fig:ch_line_results}}
\footnote{This figure was obtained from~\cite{darbon2016algorithms} without modification under the terms of the Creative Commons Attribution 4.0 International License (CC BY 4.0, \url{https://creativecommons.org/licenses/by/4.0/}).} {where the initial condition is $J(x) = \frac{1}{2}\|x\|_1^2$ and the Hamiltonian is $H(x) = \sqrt{\langle x, Dx\rangle}$ with 
certain diagonal positive definite matrix $D$. }

The Hopf representation formula is further generalized to linear dynamics, which is proposed as a future direction in~\cite{darbon2016algorithms} and then developed in~\cite{kirchner2017time}. To be specific, they consider the class of optimal control problems with linear dynamics in the following form
\begin{equation*}
\nu(y, t) = \min\left\{\int_t^T L(u(s), s) ds + J(x(T)) \colon \dot x(s) = Ax(s) + B(s) u(s) \:\forall s\in [t,T], x(t) = y\right\}. 
\end{equation*}
The function $x,t\mapsto \phi(x,t)= \nu(x, T-t)$ solves the HJ PDE with initial condition $J$ and the Hamiltonian defined by
\begin{equation*}
\begin{split}
H(x,p,t) &= \sup_{u\in \R^m} \{-\langle p, Ax + B(T-t) u \rangle - L(u, T-t) \} \\
&= -\langle p, Ax\rangle + \sup_{u\in \R^m} \{\langle -B(T-t)^T p,  u \rangle - L(u, T-t) \} \\
&= -\langle p, Ax\rangle + L^*(-B(T-t)^T p,T-t),
\end{split}
\end{equation*}
where $L^*(\cdot,t)$ is the Fenchel-Legendre transform of $L(\cdot, t)$.
To solve the HJ PDE and the corresponding optimal control problem, a change of variable $z(t) = e^{-tA} x(t)$ is applied to the optimal control formulation to get
\begin{equation*}
\nu(y, t) = \min\left\{\int_t^T L(u(s), s) ds + J(e^{TA}z(T)) \colon \dot z(s) = e^{-sA}B(s) u(s) \;\forall s\in [t,T], z(t) = e^{-tA}y\right\}. 
\end{equation*}
After setting $\tilde \nu (e^{-tA}y, t) = \nu(y,t)$, the function $x,t\mapsto \tilde \phi(x,t)= \tilde \nu(x, T-t)$ solves the HJ PDE with an initial condition $z\mapsto \tilde J(z):= J(e^{TA}z)$ and Hamiltonian 
\begin{equation*}
\tilde H(p,t) = \sup_{u\in \R^m} \left\{-\langle p, e^{-(T-t)A}B(T-t) u \rangle - L(u, T-t) \right\} = L^*(-B(T-t)^T e^{-(T-t)A^T} p, T-t).
\end{equation*}
In this case, the Hamiltonian $\tilde H$ only depends on the momentum variable $p$ and the time variable $t$. Therefore, the viscosity solution to the HJ PDE can be represented by the generalized Hopf formula:
\begin{equation*}
\tilde \phi(x,t) = -\min_{p\in \Rn}\left\{\tilde J^*(p) + \int_0^t H(p,s)ds - \langle x, p\rangle \right\},
\end{equation*}
which is solved by discretizing over time axis and then applying some optimization method, such as a relaxed Newton’s method in~\cite{kirchner2017time}. Moreover, in~\cite{kirchner2017time} the authors also studied differential games with linear dynamics, which corresponds to more complicated non-convex Hamiltonians. Following this line, in~\cite{chow2017algorithm,chow2018algorithm,kirchner2018primal}, different optimization methods, such as primal-dual method or coordinate descent method, have been applied to the generalized Hopf formula to solve more complicated non-convex Hamiltonians or differential game problems. More complicated problems have been solved with a more general form of Hopf formula, such as HJ PDEs with state-dependent non-convex Hamiltonian~\cite{chow2019algorithm}, and HJ equation in density spaces~\cite{chow2019algorithmdensity}.

Besides these general formulas, some other special cases are studied. For instance,~\cite{chen2024lax} proposes the representation formula for the class of HJ PDEs with Hamiltonian in the form of $H(p,x) = K(p) - \frac{1}{2}(x-v_0)^T M (x-v_0)$ for piecewise affine positively 1-homogeneous convex function $K$, symmetric positive definite matrix $M$, and a vector $v_0\in \Rn$, with assumptions relating $M$ to the level set of $K$. These HJ PDEs correspond to optimal control problems with a quadratic running cost on the state variable and a rectangular constraint on velocity, which makes the optimal trajectory in a piecewise linear form. Then, the viscosity solution can be represented by a minimization problem of a piecewise third order polynomial function adding with the initial condition. The dual case is presented in~\cite{chen2024hopf}, which considers the HJ PDEs with Hamiltonian $H(p,x) = \frac{1}{2}\langle p, Mp\rangle - U(x)$, where $U$ is a positive 1-homogeneous convex function whose level set is determined by the matrix $M$. This corresponds to the optimal control problem with running cost to be a combination of kinetic energy and potential energy $U$. Under these assumptions, the optimal trajectory is a piecewise quadratic function, and the viscosity solution to the HJ PDE is a supremum of difference of a piecewise polynomial function with the Fenchel-Legendre transform of the initial condition.

In addition to the representation formulas above, there is a technique for HJ PDEs called max-plus or min-plus algebra~\cite{mceneaney2006max}. Specifically, under certain assumptions, if the initial condition $J$ for an HJ PDE is the minimum of several functions $J_1, \dots, J_M$, i.e., $J = \min_{i=1,\dots, M} J_i$, then the solution $S$ to the HJ PDE with initial condition $J$ is also the minimum of the corresponding solutions $S_1, \dots, S_M$ to the HJ PDEs with initial conditions $J_1, \dots, J_M$, i.e., we have $S = \min_{i=1,\dots, M} S_i$. This property has been used in~\cite{darbon2016algorithms,chen2024lax,chen2024hopf,darbon2023neural,gaubert2011curse,fleming2000max,dower2015max} to solve HJ PDEs with more complicated initial conditions.

\section{Monte-Carlo Sampling via Laplace's Method}
\label{sec: Laplace_approximations}
Laplace's method is a classical analytical technique often used in asymptotic analysis to approximate integrals~\cite{laplace1774memoire}. 
In the context of HJ PDEs, this method can be employed to derive approximate solutions in cases where the solution exhibits an exponential form, commonly seen in problems of optimal control or short-time asymptotics.
Formally, Laplace's method can be explained as follows. Let $\varphi \colon \mathbb{R}^d \to \mathbb{R}$ be twice continuously differentiable and $h \colon \mathbb{R}^d \to \mathbb{R}$ be a continuous function. Laplace's method provides an asymptotic equivalence for an integral that becomes increasingly
peaked around the global minimizer $u^\star$ of $\varphi$, assumed to be unique:
\begin{equation}
    \int h(u) \exp(-\varphi(u)/\delta)du \sim \frac{(2 \pi \delta)^d}{\exp(\varphi(u^\star)/\delta)} \frac{h(u^\star)}{\sqrt{\det(\nabla^2 \varphi(u^\star))}},
    \label{eq: laplace_method}
\end{equation}
where $a(\delta) \sim b(\delta)$ is a short hand for $a(\delta)/b(\delta) \to 1$ as $\delta \to 0$.
The primary implication of~\eqref{eq: laplace_method} is that one can use Laplace's method in two different ways: one may estimate the global minimizer of a function $\varphi$ by estimating an integral or vice-versa. The latter is a common technique found in statistical inference~\cite{kaipio2006statistical}. In most situations, either approach can be difficult, as one must either estimate a high-dimensional integral or solve a non-convex optimization problem.

Recent works have studied the connection between Laplace's method and infimal convolutions~\cite{tibshirani2024laplace} with applications in zeroth-order~\cite{osher2023hamilton, di2025monte} and global~\cite{heaton2024global, zhang2024inexact} optimization. 
It also provides an approach (via self-normalization) to avoid computing the determinant seen in~\eqref{eq: laplace_method}. In particular, \emph{one may estimate the minimizer of $\varphi$ as} 
\begin{equation}
    u^\star \sim \dfrac{\int u \exp(-\varphi(u)/\delta) du}{\int \exp(-\varphi(u)/\delta) du}.
    \label{eq: laplace_normalized}
\end{equation}

{The intuition behind this why this expression leads to a minimizer is that the weight $\exp(-\varphi(u)/\delta)$ is sharply peaked around the minimizer of $\varphi$. 
As $\delta \to 0$, the probability mass of this distribution concentrates near $u^\star$, so the normalized integral acts like a \emph{weighted average around the peak}, collapsing to the true minimizer in the small-$\delta$ limit. 
In this way,~\eqref{eq: laplace_normalized} connects optimization with sampling: instead of solving $\nabla \varphi(u) = 0$ directly, one can estimate $u^\star$ by drawing Monte Carlo samples from the density
\begin{equation}
    p_\delta(u) \propto \exp(-\varphi(u)/\delta).
\end{equation}
}
{
As an explicit (but simple) example, if $\varphi(u) = (u-\mu)^2$ on $\mathbb{R}$, then $p_\delta(u) \propto \exp(-(u-\mu)^2/\delta)$ is exactly Gaussian with mean $\mu$ and variance $\delta/2$. In this case, the formula~\eqref{eq: laplace_normalized} gives
\begin{equation}
    \frac{\int u \exp(-(u-\mu)^2/\delta) \, du}{\int \exp(-(u-\mu)^2/\delta) \, du} = \mu,
\end{equation}
which coincides precisely with the true minimizer of $\varphi$. 
Thus, for quadratic $\varphi$, the Laplace estimator is not just asymptotic---it is exact for all $\delta > 0$.
}

In the context of HJ PDEs, this normalized version of Laplace's method can be used to approximate the solution to the Lax-Oleinik formula given in~\eqref{eq: lax_oleinik} by setting $\varphi(u) = J(u) + tH^\star\left(\frac{x-u}{t}\right)$.
After approximating $u^\star$ via Laplace's using~\eqref{eq: laplace_normalized}, we may plug it back into the objective function to obtain $S(x,t)$.~\cite{osher2023hamilton} shows convergence of~\eqref{eq: laplace_normalized} in the setting where $J(u)$ is $t$-weakly convex with $\frac{1}{t} H^\star(\frac{x-u}{t}) = \frac{1}{2t}\|x-u\|^2$. This result is generalized in~\cite[Theorem 2]{tibshirani2024laplace} to the case where $\varphi$ is Hölder continuous.

\begin{figure}
    \centering
    \begin{tabular}{ccc}
        \includegraphics[width=0.27\textwidth, height=1.1in]{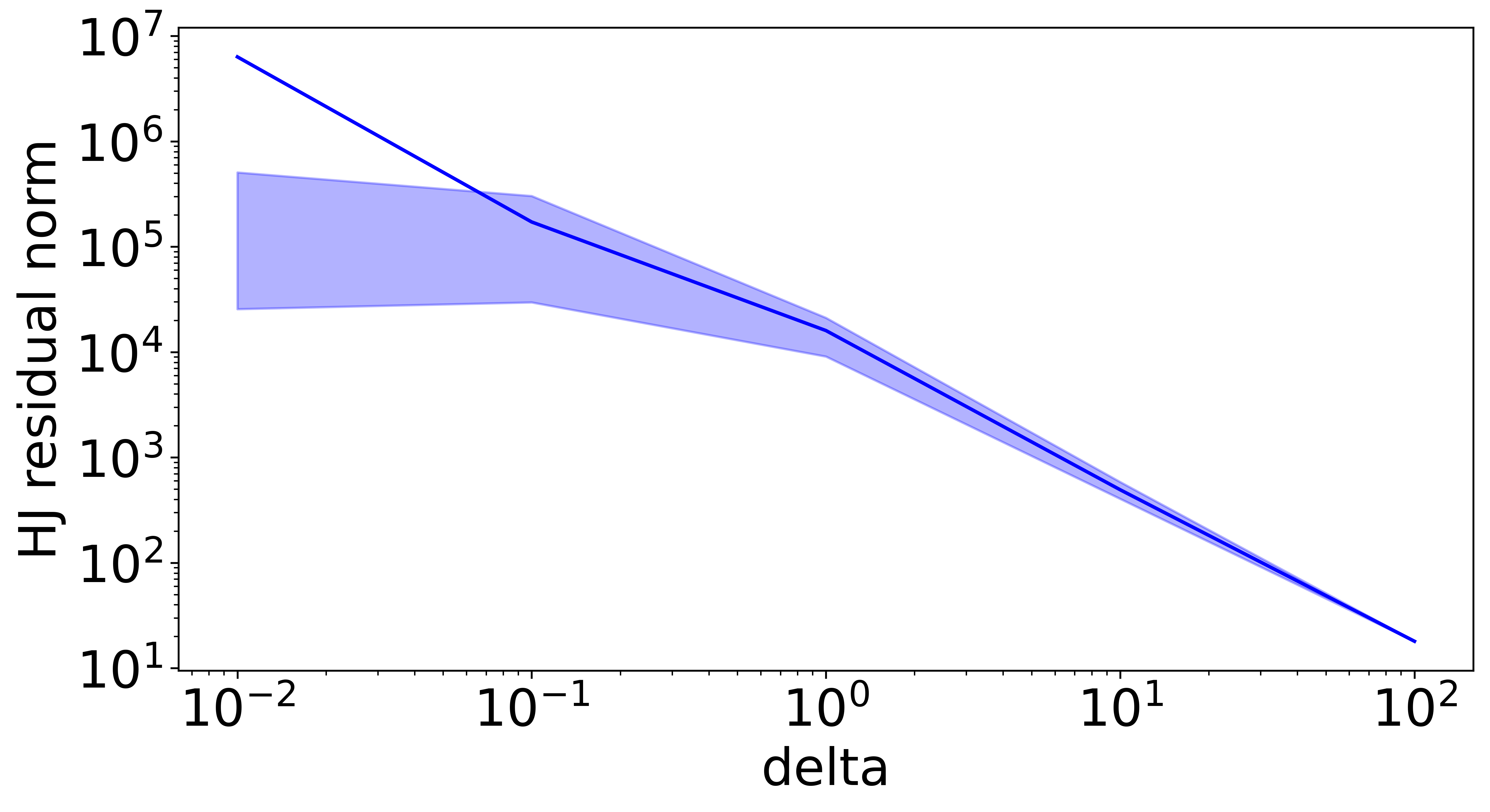}
        & 
        \includegraphics[width=0.27\textwidth, height=1.1in]{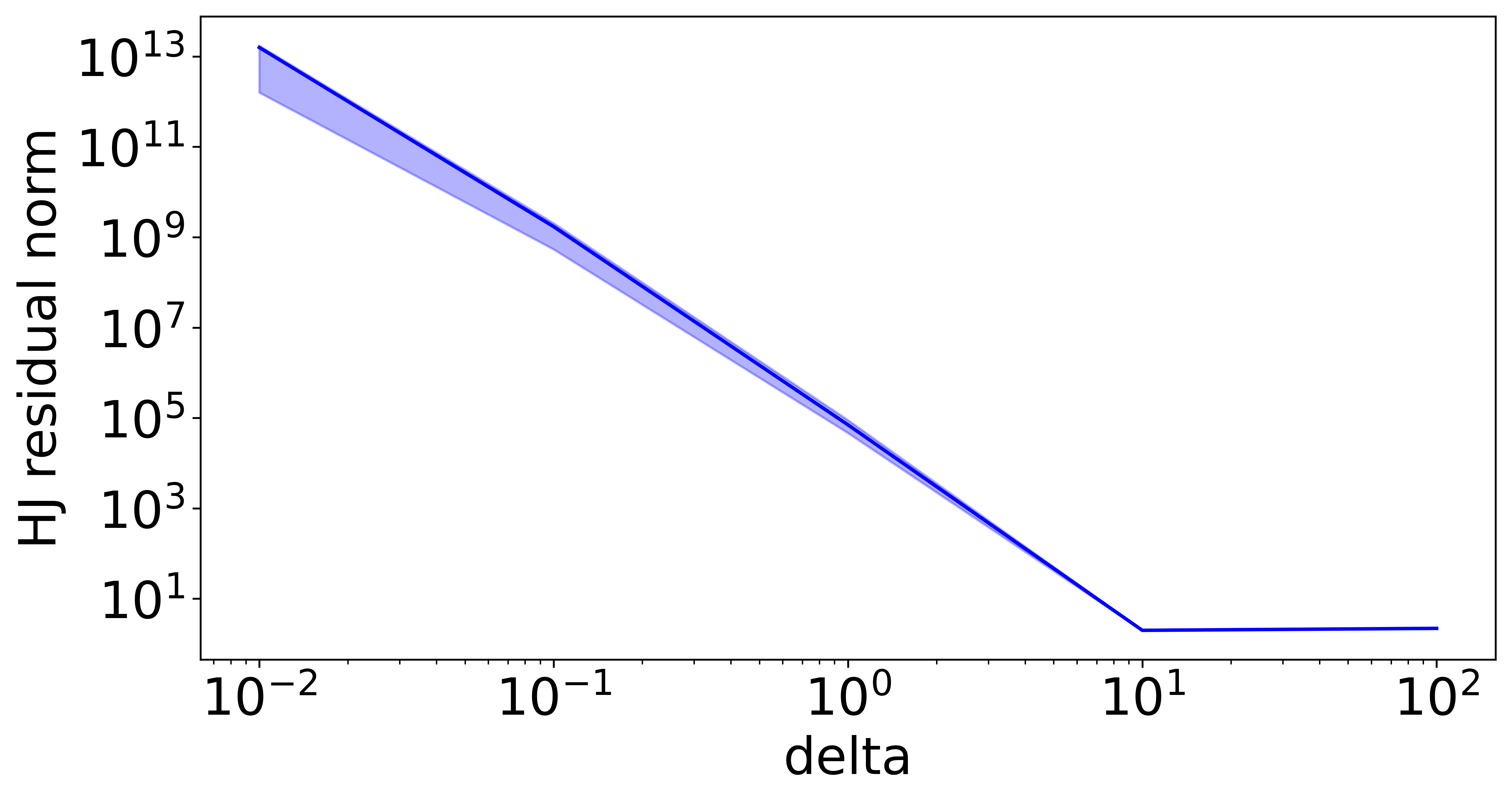}
        &
        \includegraphics[width=0.27\textwidth, height=1.1in]{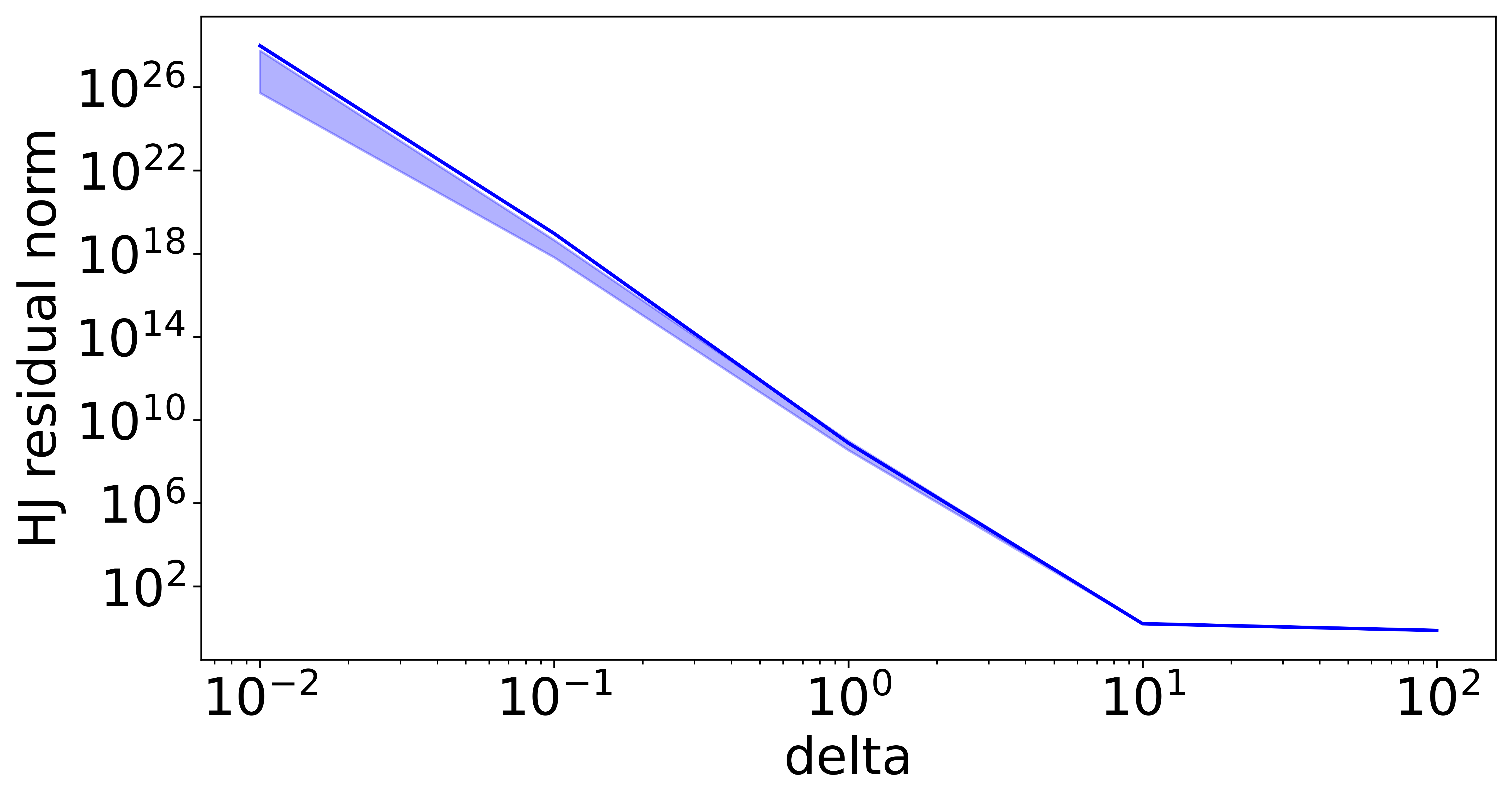}
    \end{tabular}
    \caption{\small{Median along with the 20th and 80th percentiles of HJ residuals in dimension 10 and Hamiltonians given by $H = \frac{1}{p} \| \cdot \|_p^p$, p = {1, 5, 10}, from left to right.} In these experiments, $10^5$ samples are used for the Monte Carlo approximation of the integrals in~\eqref{eq: laplace_normalized}.}
    \label{fig: laplace_HJ_experiments}
\end{figure}

Thus, Laplace's method offers a different perspective based on Monte Carlo techniques for approximating solutions to HJ PDEs, setting it apart from grid-based discretization or optimization of the Hopf and Lax-Oleinik formulas. 
Unlike traditional numerical methods that rely on fixed grids and often struggle with scalability in high-dimensional settings, Laplace's method avoids the use of grids entirely by employing probabilistic sampling to approximate the integral~\eqref{eq: laplace_normalized}. 
In particular, Laplace's method can potentially be more efficient if approximation of the integral is easier to compute than traditional methods, e.g., ENO and WENO~\cite{osher1988fronts}. Results based on~\cite{tibshirani2024laplace} are shown in Fig.~\ref{fig: laplace_HJ_experiments}, where the effect of the choice of the smoothing parameter $\delta$ is varied for HJ equations with $H(\cdot) = \frac{1}{p}\|\cdot\|_p^p$ for $p = 2, 5,$ and $10$.
We note these connections between infimal convolutions and HJ PDEs have also been explored in the context of posterior mean estimation and denoising~\cite{darbon2021connecting, darbon2021bayesian}.

\section{Deep Learning Methods}
\label{sec: deep_learning}
Deep learning methods have gained traction in solving high-dimensional HJ PDEs. Neural networks serve as flexible function approximators, leveraging optimization frameworks to approximate solutions. Without the need for discretization, these neural network methods are able to mitigate the curse of dimensionality.
While deep learning methods for solving general PDEs~\cite{beck2023overview} can be used here, we focus on works that are more specific to HJ PDEs. \update{Beyond the approaches we will survey in this section, other machine learning methods have also been proposed for related PDE and control problems (see, e.g.,~\cite{bansal2021deepreach,esteve2025finite,kim2021hamilton}).}

\update{Unlike classical machine learning where methods are often categorized as supervised 
versus unsupervised depending on the availability of labeled data, most scientific 
machine learning (SciML) approaches for HJ PDEs do not fall neatly into this dichotomy. 
Rather than relying on external datasets, these methods derive their training signals 
from the governing equations themselves. In what follows, we highlight the type of 
signal each method relies on: Section~\ref{sec:ML_AC} focuses on loss 
functions motivated by reinforcement-learning style objectives, while 
Section~\ref{sec:ML_representation} emphasizes neural network architectures informed by 
representation formulas. Moreover, all methods surveyed in this section have published 
code, with links provided in their original papers.}

\subsection{Actor-Critic Methods}\label{sec:ML_AC}
\update{Actor–critic algorithms are a standard class of methods in deep reinforcement learning. 
They combine a critic, which estimates the value function of a given policy, typically 
through temporal-difference (TD) learning, with an actor that updates the policy in the 
direction suggested by the critic’s feedback. This canonical two-step template has been 
adapted in various ways to scientific computing problems.}

\update{Building on this framework,} in \cite{zhou2021actor}, the authors propose a numerical method for solving high-dimensional static HJ equations using an actor-critic framework combined with neural networks. \update{The code of this method is released in~\url{https://github.com/MoZhou1995/DeepPDE\_ActorCritic}.} The convergence of this method is analyzed in~\cite{zhou2023policy,zhou2024solving}. The central equation they tackle takes the form:

\begin{equation*}
    \inf_{u \in U} \left[\frac{1}{2} \text{Tr}(\update{\sigma(x,u)\sigma(x,u)^\top \text{Hess}(V)(x)}) + b(x,u)^\top\nabla V(x) + f(x,u)\right] - \gamma V(x) = 0 \quad \text{in } \Omega,
\end{equation*}
where $V$ is the value function, $u$ is the control, and $\gamma$ is the discount rate. Their key methodological contribution is a variance-reduced least squares temporal difference (VR-LSTD) method for policy evaluation. The critic step involves minimizing a loss functional $L_1(V) = \mathbb{E}[\text{TD}_1^u]^2$ with a quadratic loss term on boundary condition, where
\begin{equation*}
\begin{adjustbox}{width=0.99\textwidth}$
    \text{TD}_1^u = \int_0^T e^{-\gamma s}f(X_s,u(X_s))ds - \int_0^T e^{-\gamma s}\nabla V(X_s)^\top\sigma(X_s,u(X_s))dW_s + e^{-\gamma T}V(X_T) - V(X_0).
$\end{adjustbox}
\end{equation*}
The actor step involves the following loss:
\begin{equation*}
\mathbb{E}_{X_0 \sim \mu, u} \left[ \int_{0}^{T \wedge \tau} f(X_s, u(X_s)) e^{-\gamma s} \,\mathrm{d}s 
+ V(X_{T \wedge \tau}) e^{-\gamma (T \wedge \tau)} \right],
\end{equation*}
\update{where $\mu$ is some initial distribution for $X_0$ and $\tau$ is the stopping time when the process hits the domain boundary.}
They show that this VR-LSTD method has better convergence properties than standard LSTD approaches.

\begin{figure*}[t!]
    \centering
    \includegraphics[width=0.98\textwidth]{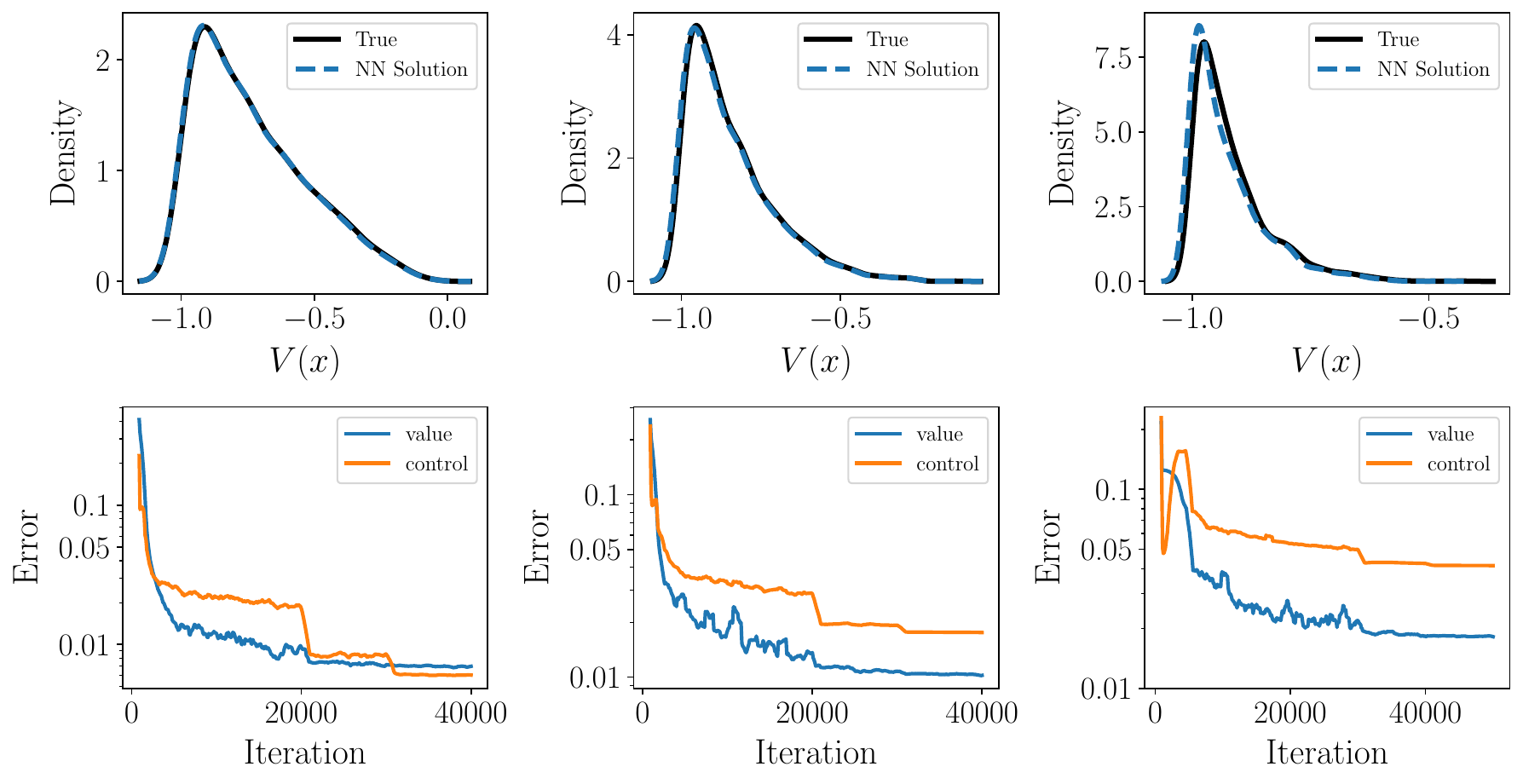}
    \caption{Top: density of the value of $V$ for the \update{diffusive Eikonal equation~\eqref{eqt:ACmethod_Eikonal} (with $a_2=1.2$, $a_3=0.2$, and $R=1$)} with $d=5$ (left), $d=10$ (middle), and $d=20$ (right).
    Bottom: associated error curves in the training process with $d=5$ (left), $d=10$ (middle), and $d=20$ (right).}
    \label{fig:ekn}
\end{figure*}

To improve numerical accuracy near domain boundaries, they introduce an adaptive step size scheme when the sampled state is near the boundary.
This scheme, combined with their neural network parameterization of both value and control functions, enables them to solve the equation up to 20 dimensions with relative errors typically around 1-5\%. They demonstrate the effectiveness of their approach on several challenging problems including linear quadratic regulators, stochastic Van der Pol oscillators, and diffusive Eikonal equations, showing significant improvement over traditional numerical methods that struggle with the curse of dimensionality.
For instance, they tested on the following diffusive Eikonal equation:
\begin{equation}\label{eqt:ACmethod_Eikonal}
\left\{ \begin{aligned}
& \epsilon \Delta V(x) + \inf_{u \in B_1} \bigl(c(x) u^{\top} \nabla V(x)\bigr) + 1  = 0 ~~~ \text{in} ~ B_R, \\
&V(x) = a_3-a_2 ~~~ \text{on} ~ \partial B_R,
\end{aligned} \right.
\end{equation}
where $c(x)$ is set to be $c(x) = \dfrac{3(d+1)a_3}{2d a_2(2a_2 - 3a_3|x|)} > 0$. The numerical results are shown in Fig.~\ref{fig:ekn}. This methodology also applies to other formulations of control problems, such as finite-horizon settings, time-inhomogeneous systems, control without discounting, and problems without boundary conditions.

Another actor-critic approach~\cite{lin2021alternating} solves HJ PDEs with a viscosity term on the right-hand side. Here, a neural network method for solving high-dimensional stochastic mean-field games (MFGs) is proposed. Their framework tackles the coupled system of HJ PDEs and Fokker-Planck (FP) equations that characterize MFGs:
\begin{equation}\label{eqt:HJ_FP}
    \begin{split}
    - \partial_t \phi - \nu \Delta \phi + H(x, \nabla \phi) = f(x, \rho), \quad\quad
    \partial_t \rho - \nu \Delta \rho - \text{div}(\rho \nabla_p H(x, \nabla \phi)) = 0,
    \end{split}
\end{equation}
with initial condition $\rho(x,0) = \rho_0(x)$ and terminal condition $\phi(x,T) = g(x,\rho(\cdot,T))$. 
The authors reformulate this system as a convex-concave saddle-point problem, which they express as:
\begin{equation*}
    \begin{split}
    \inf_{\rho(x,0)=\rho_0(x)} \sup_{\phi} \int_0^T \int_{\Omega} (\partial_t \phi + \nu \Delta \phi - H(x, \nabla \phi))\rho(x,t) \, dx + F(\rho(\cdot,t)) \, dt 
    \\
    + \int_{\Omega} \phi(x,0)\rho_0(x) \, dx - \int_{\Omega} \phi(x,T)\rho(x,T) \, dx + G(\rho(\cdot,T)),
    \end{split}
\end{equation*}
where the first variations of $F$ and $G$ are the functions $f$ and $g$ in~\eqref{eqt:HJ_FP}.
This variational structure allows them to parameterize $\rho$ and $\phi$ (the actor and the critic) with neural networks and train them in an alternating fashion, similar to generative adversarial networks (GANs)~\cite{goodfellow2020generative}. The resulting algorithm can approximate solutions to high-dimensional HJ PDEs, and in their experiments, HJ PDEs of dimensions up to 100 are solved efficiently~\cite[Section 5]{lin2021alternating}.

\update{More broadly, it is worth noting that essentially all value-based reinforcement 
learning algorithms can be interpreted as numerical schemes for HJ-type equations. 
Examples include deep Q-learning method~\cite{kim2021hamilton}, 
which has also been applied to control problems closely connected with HJ equations.
}

\subsection{Neural Network Architectures Inspired by Representation Formulas}\label{sec:ML_representation}
In this section, we review several works that propose neural network architectures inspired by representation formulas for HJ PDEs under different assumptions. These architectures, when assigned appropriate parameters based on their theoretical correspondence, inherit rigorous guarantees from HJ theory—eliminating the need for a traditional training process. 


In~\cite{darbon2020overcoming,darbon2021some,darbon2023neural}, different neural network architectures are designed based on the representation formulas of solution to certain classes of HJ 
PDEs. 
In~\cite{darbon2020overcoming}, inspired by the Hopf formula~\eqref{eqt:Hopf}, a shallow neural network with a max activation function is proposed to solve HJ PDEs. Specifically, the formula for the neural network is
\[
\phi(x, t; \{(p_i, \theta_i, \gamma_i)\}_{i=1}^{m}) = \max_{i =1, \dots, m} \left\{ \langle p_i, x \rangle - t\theta_i - \gamma_i \right\},
\]
where $\{(p_i, \theta_i, \gamma_i)\}_{i=1}^{m}$ denotes the parameters in the neural network. Moreover, it is proven that, under certain technical assumptions, if there exists a convex interpolation for $(p_i, \gamma_i)$ and the points $p_1,\dots, p_m$ are distinct, then the function represented by the neural network is the viscosity solution to the HJ PDE with the initial condition $J$ given by
\[
J(x) = \max_{i =1, \dots, m} \left\{ \langle p_i, x \rangle - \gamma_i \right\},
\]
and the Hamiltonian $H$ defined as 
\begin{equation}\label{eqt:NN_Hamiltonian}
H({p}) =
\begin{dcases}
\inf\limits_{\alpha \in \mathcal{A}({p})} \left\{ \sum_{i=1}^{m} \alpha_i \theta_i \right\}, & \text{if } {p} \in \operatorname{dom} J^*, \\
+\infty, & \text{otherwise}.
\end{dcases}
\end{equation}
Here, the set $\mathcal{A}({p})$ is characterized as 
\[
\mathcal{A}({p}) = \argmin_{\substack{\alpha_1, \dots, \alpha_m\geq 0\\
\sum_{i=1}^m \alpha_i = 1\\ \sum_{i=1}^{m} \alpha_i {p}_i = {p}}} \left\{ \sum_{i=1}^{m} \alpha_i \gamma_i \right\}.
\]
Moreover, for a fixed set of parameters \(\{(p_i, \theta_i, \gamma_i)\}_{i=1}^{m}\), the neural network representation is capable of solving multiple HJ PDEs. Specifically, it solves the HJ PDE with a spatially and temporally independent Hamiltonian \(\tilde H\) if and only if the conditions
$\tilde{H}(p_i) = H(p_i) \, \forall i = 1, \dots, m$
and $\tilde{H}(p) \geq H(p) \, \forall p \in \operatorname{dom} J^*$
are both satisfied. Here, the function $H$ is the Hamiltonian defined in~\eqref{eqt:NN_Hamiltonian}.


In~\cite{darbon2021some}, based on the Lax-Oleinik formula~\eqref{eq: lax_oleinik}, the authors propose two neural network architectures that solve two classes of HJ PDEs, respectively. The first architecture is designed as 
\begin{equation*}
\phi({x}, t; \{(u_i, a_i)\}_{i=1}^m) = \min_{i =1, \dots, m} \left\{ t L \left( \frac{x - u_i}{t} \right) + a_i \right\}.
\end{equation*}
This is a two-layer neural network employing an activation function \( L \) and a minimum operator. Under certain assumptions, this architecture computes the viscosity solution of an HJ PDE whose Hamiltonian is given by the convex conjugate \( L^* \), independent of the spatial and temporal variables \((x,t)\). The corresponding initial condition is expressed as the minimum of shifts of the asymptotic function associated with \( L \). Similarly, the second architecture utilizes an activation function \( J \) combined with a minimum operator. Under suitable conditions, this network solves the HJ PDE characterized by a piecewise affine Hamiltonian and a concave initial condition given by \( J \).

\begin{figure}[htbp]
\begin{minipage}[b]{.24\linewidth}
  \centering
  \includegraphics[width=.99\textwidth]{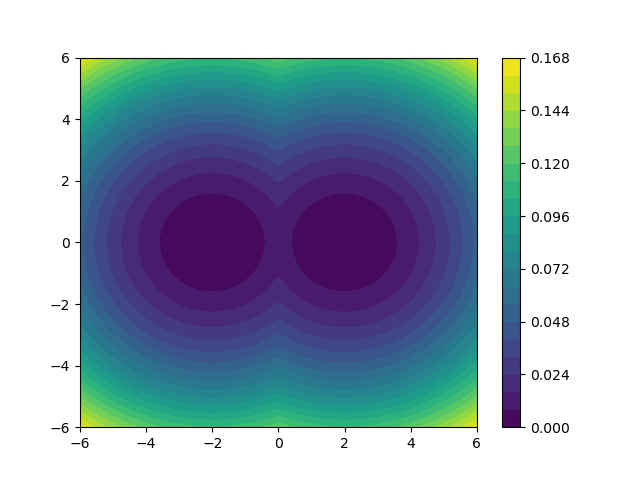}
  \centerline{\footnotesize{(a) $t=1$}}\medskip
\end{minipage}
\begin{minipage}[b]{.24\linewidth}
  \centering
  \includegraphics[width=.99\textwidth]{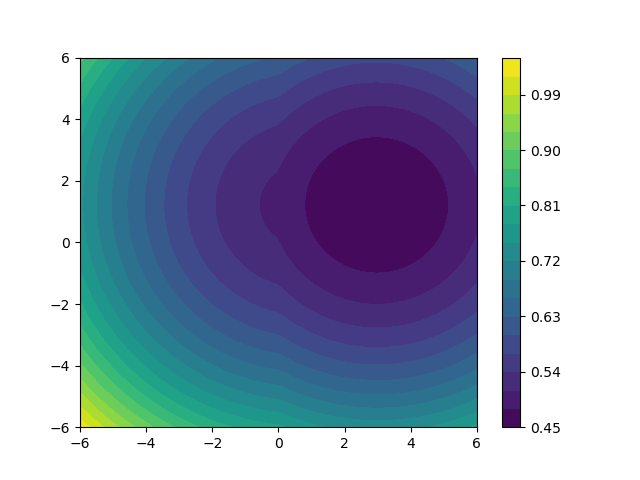}
  \centerline{\footnotesize{(b) $t=0.995$}}\medskip
\end{minipage}
\begin{minipage}[b]{.24\linewidth}
  \centering
    \includegraphics[width=.99\textwidth]{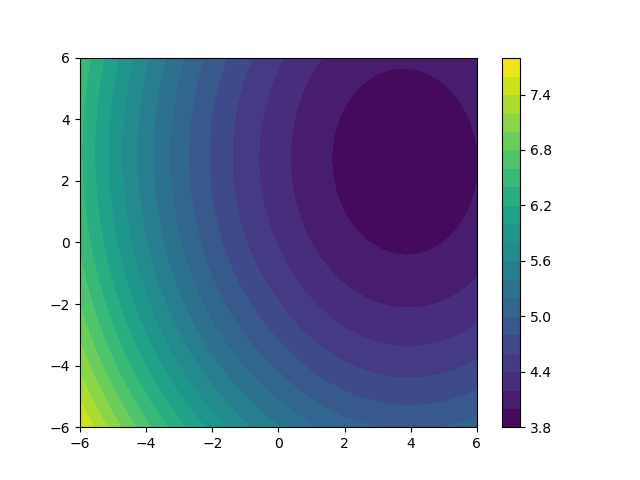}
  \centerline{\footnotesize{(c) $t=0.95$}}\medskip
\end{minipage}
\begin{minipage}[b]{.24\linewidth}
  \centering
    \includegraphics[width=.99\textwidth]{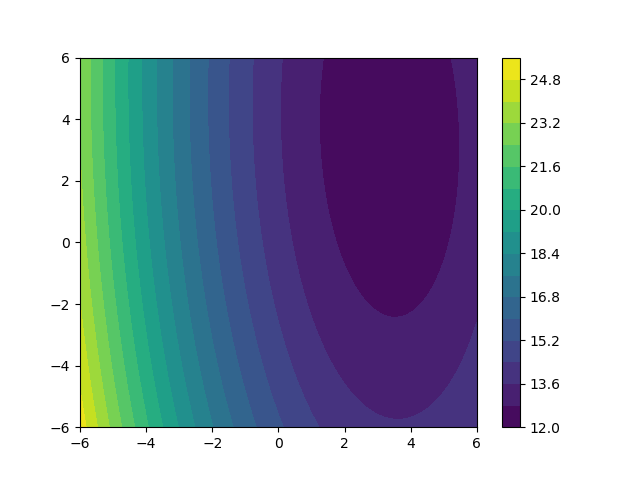}
  \centerline{\footnotesize{(d) $t=0.75$}}\medskip
\end{minipage}
\caption{
Visualization of the viscosity solution to the 16-dimensional HJ PDE with initial condition~\eqref{eqt: newton_J}, corresponding to the running cost \(L\) defined by~\eqref{eqt:eg3_L}, dynamics \(f\) defined by~\eqref{eqt:eg3_f}, and terminal time \(T = 1\).  
Two-dimensional slices at times \( t=1 \) (terminal cost), \( t=0.995 \), \( t=0.95 \), and \( t=0.75 \) are displayed in subfigures (a), (b), (c), and (d), respectively.  
Each subfigure shows solution values (indicated by color) at points of the form \((x_1,\mathbf{0}, x_2,\mathbf{0})\in\mathbb{R}^{16}\), where \(x_1\) and \(x_2\) correspond to the horizontal and vertical axes, and \(\mathbf{0}\) denotes a 7-dimensional zero vector.
\label{fig:nn_minplus_newton}}
\end{figure}
\begin{figure}
    \begin{minipage}[b]{.24\linewidth}
  \centering
    \includegraphics[width=.99\textwidth]{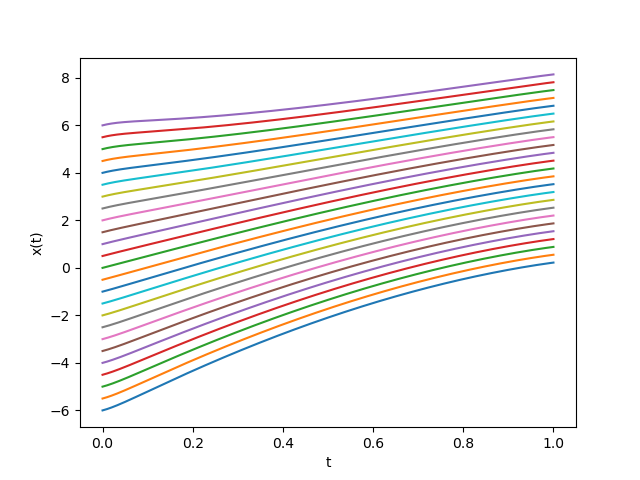}
  \centerline{\footnotesize{(a) $x_1^*$, $T=1$}}\medskip
\end{minipage}
\begin{minipage}[b]{.24\linewidth}
  \centering
    \includegraphics[width=.99\textwidth]{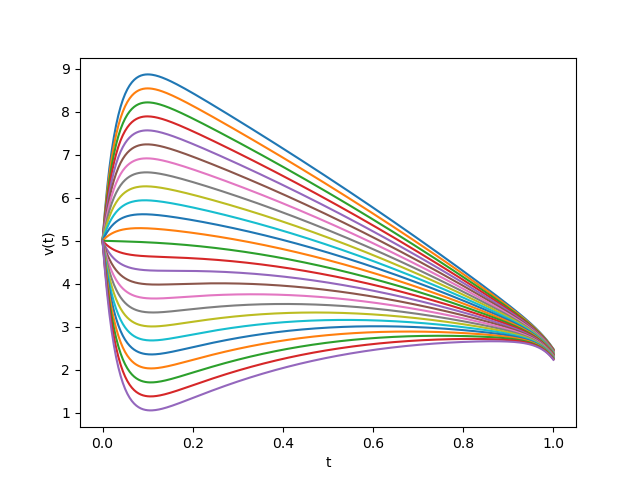}
  \centerline{\footnotesize{(b) $x_9^*$, $T=1$}}\medskip
\end{minipage}
\begin{minipage}[b]{.24\linewidth}
  \centering
    \includegraphics[width=.99\textwidth]{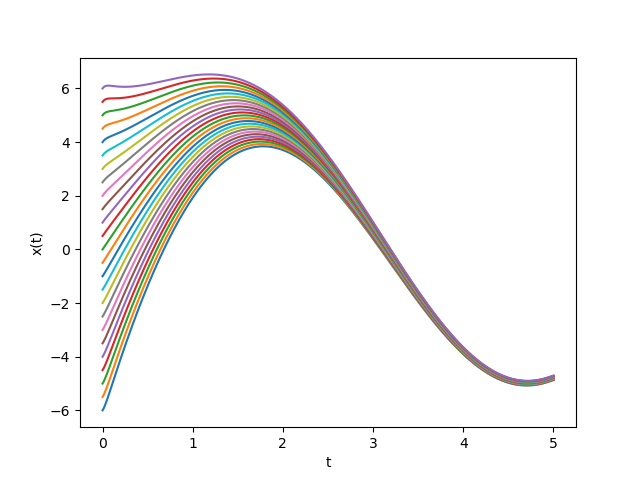}
  \centerline{\footnotesize{(c) $x_1^*$, $T=5$}}\medskip
\end{minipage}
\begin{minipage}[b]{.24\linewidth}
  \centering
    \includegraphics[width=.99\textwidth]{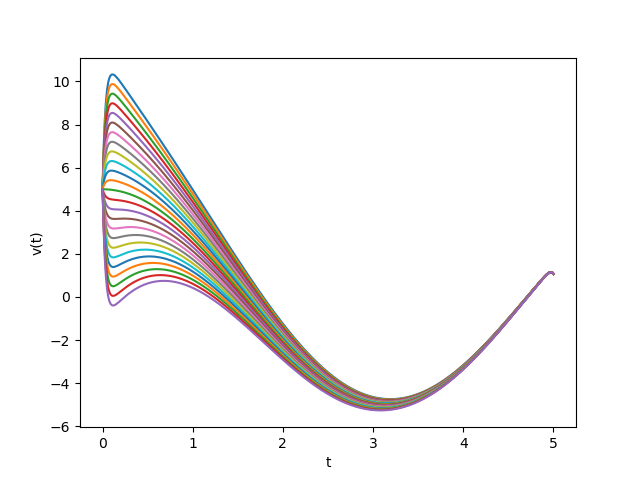}
  \centerline{\footnotesize{(d) $x_9^*$, $T=5$}}\medskip
\end{minipage}
\begin{minipage}[b]{.24\linewidth}
  \centering
    \includegraphics[width=.99\textwidth]{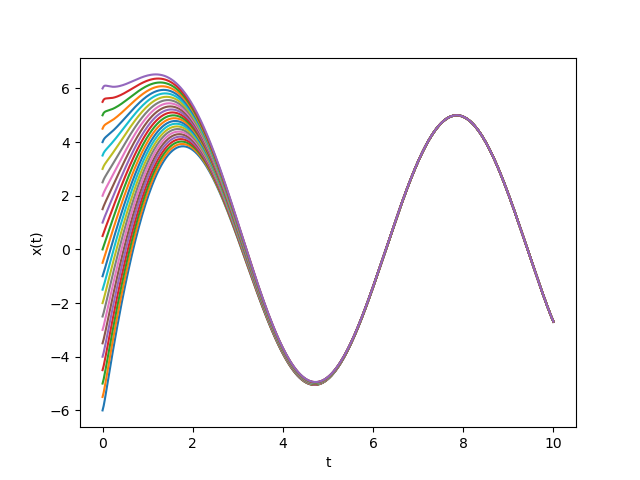}
  \centerline{\footnotesize{(e) $x_1^*$, $T=10$}}\medskip
\end{minipage}
\begin{minipage}[b]{.24\linewidth}
  \centering
    \includegraphics[width=.99\textwidth]{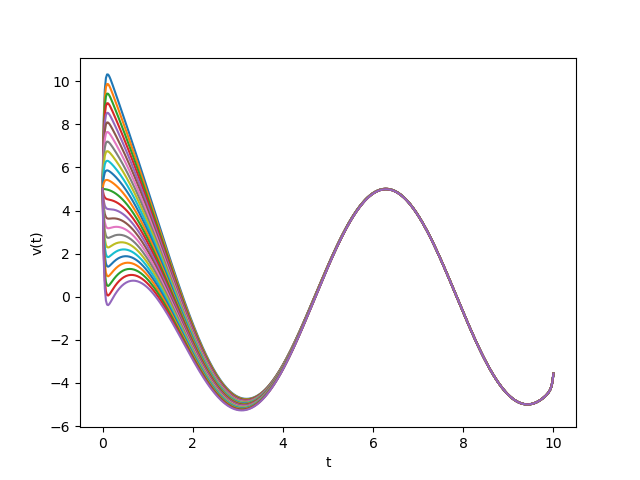}
  \centerline{\footnotesize{(f) $x_9^*$, $T=10$}}\medskip
\end{minipage}
    \caption{
    Visualization of optimal trajectories \( x^* \) for the 16-dimensional optimal control problem with running cost \( L \) given by~\eqref{eqt:eg3_L}, dynamics given by~\eqref{eqt:eg3_f}, terminal cost~\eqref{eqt: newton_J}, and terminal times \(T = 1,\,5,\,10\). Panels (a), (c), and (e) display the evolution of the first component \( x_1^* \) of the optimal trajectories, while panels (b), (d), and (f) display the corresponding evolution of the ninth component \( x_9^* \). In each panel, different curves represent trajectories originating from different initial conditions \( x_0 = (x,\mathbf{0}) \in \mathbb{R}^{16} \), where \(\mathbf{0}\) is the zero vector in \(\mathbb{R}^{15}\).
    \label{fig:nn_minplus_xv}}
\end{figure}

In~\cite{darbon2023neural}, the authors propose a ResNet-type deep neural network based on the linear-quadratic regulator framework and min-plus algebra to solve HJ PDEs. These PDEs have quadratic Hamiltonians $H(p,x,t)$ with respect to $(p,x)$, whose coefficients depend on time $t$, and initial conditions expressed as the minimum of several quadratic functions. Additionally, the authors design a neural network architecture to compute the optimal control and corresponding optimal trajectories associated with these control problems.
The numerical results for HJ PDEs exhibiting Newtonian dynamics, originally presented in~\cite{darbon2023neural}, are reproduced in Figs.~\ref{fig:nn_minplus_newton} and~\ref{fig:nn_minplus_xv}.
Specifically, the simulations consider the dimension setting $n=2m=16$, with the Lagrangian $L$ defined by
%
%
\begin{equation}\label{eqt:eg3_L}
    L(u,x,t) = \frac{1}{2} \left\|x - 5\sin t \begin{pmatrix}
        \mathbf{1} \\
        \mathbf{0}
    \end{pmatrix} + 5\cos t \begin{pmatrix}
        \mathbf{0} \\
        \mathbf{1}
    \end{pmatrix}\right\|^2 + \frac{1}{2000} \|u\|^2,
\end{equation}
and the source term $f$ given by
\begin{equation}\label{eqt:eg3_f}
    f(u,x,t) = \begin{pmatrix}
        O & I \\
        O & O
    \end{pmatrix}x+ \begin{pmatrix}
        O \\
        I
    \end{pmatrix}u,
\end{equation}
for $u\in \R^m$, $x\in \Rn$, and $t\in [0,T]$.
Here, $\mathbf{1}$ and $\mathbf{0}$ represent the one and zero constant vectors in $\mathbb{R}^m$, while $O$ and $I$ denote the $m\times m$ zero and identity matrices, respectively.
The terminal cost function $J\colon \mathbb{R}^n\to\mathbb{R}$ is defined as
\begin{equation} \label{eqt: newton_J}
    J(x) = \min\left\{\frac{1}{320}\left((x_1 + 2)^2 + \sum_{i=2}^nx_i^2\right), \frac{1}{320}\left((x_1 - 2)^2 + \sum_{i=2}^nx_i^2\right)\right\},
\end{equation}
for $x=(x_1,x_2,\dots, x_n)\in\Rn$.
Fig.~\ref{fig:nn_minplus_newton} shows two-dimensional slices of the computed HJ PDE solutions at different times. Furthermore, Fig.~\ref{fig:nn_minplus_xv} displays the optimal trajectories corresponding to terminal times $T=1$, $T=5$, and $T=10$ from various initial positions $x_0=(x,\mathbf{0})\in\mathbb{R}^{16}$ (where $\mathbf{0}\in\mathbb{R}^{15}$ is the zero vector). 


\section{Saddle Point Methods}
\label{sec: saddle_point}
In this section, we review recent numerical methods~\cite{meng2023primal,meng2024primal} for solving HJ PDEs based on saddle-point formulations. These approaches exploit connections between HJ PDEs and problems in optimal transport and mean-field games by reformulating the PDEs as components of certain saddle-point problems. 
\update{This 
reformulation is attractive because it introduces flexibility beyond standard PDE 
discretizations, allowing one to incorporate primal–dual methods, preconditioning, and 
other optimization techniques that may be advantageous in high-dimensional or 
structured settings. In principle, the resulting saddle-point problems can be approached 
in different ways; in the works cited here, they are solved after spatial discretization 
using primal–dual algorithms such as the Primal–Dual Hybrid Gradient (PDHG) 
algorithm~\cite{chambolle2011first}.}

Specifically, we consider the HJ PDE~\eqref{eqt:HJ_general} defined on a rectangular spatial domain \(\Omega = \prod_{i=1}^n [a_i, b_i]\) with periodic boundary conditions, \update{and denote the solution by $\phi$}.

\begin{figure}[htbp]
    \centering
    \begin{subfigure}{0.3\textwidth}
        \centering \includegraphics[width=\textwidth]{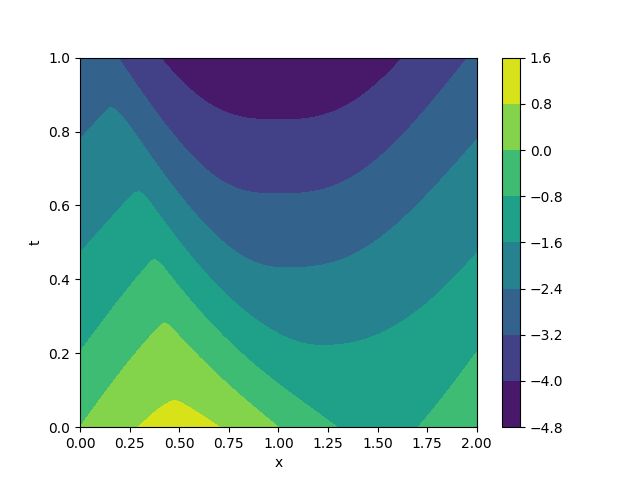}
    \end{subfigure}
    \begin{subfigure}{0.3\textwidth}
        \centering \includegraphics[width=\textwidth]{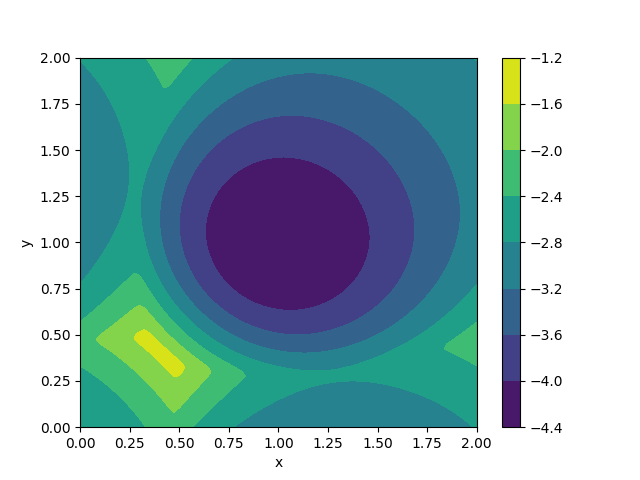}
    \end{subfigure}
    \caption{Contour plots illustrating the solution of the HJ PDE with Hamiltonian \( H(p, x) = \|p\| + 3\exp(-4\|x - 1\|^2) + 1 \) and initial condition \( J(x) = \sum_{i=1}^n \sin (\pi x_i) \).
    The left panel shows the solution for the one-dimensional case, while the right panel presents the level sets at $t=0.75$ for the two-dimensional case.
    }
    \label{fig:saddle_pt_results}
\end{figure}

In~\cite{meng2023primal}, the solution $\phi$ is represented as a part of the saddle point to the following problem
\begin{equation*}
\begin{split}
\min_{\substack{\phi\\ \phi(x,0)=J(x)}}\max_{\rho} \Bigg\{ \int_0^T \int_\Omega \rho(x,t)\left(\frac{\partial \phi(x,t)}{\partial t}  +\max_{v\in\Rn} \left\{\langle v,\nabla_x\phi(x,t)\rangle - H^*\left(v,x, t\right)\right\} \right)dxdt\\
- c\int_\Omega \phi(x,T)dx\Bigg\},
\end{split}
\end{equation*}
where $c > 0$ is a hyperparameter, and $H^*(\cdot,x,t)$ is the Fenchel-Legendre transform of $H(\cdot,x,t)$. 
If we further assume the function $\rho$ in the saddle point is non-negative, then the problem becomes
\begin{equation*}
\begin{split}
    \min_{\substack{\phi\\ \phi(x,0)=J(x)}}\max_{\substack{\rho,v\\ \rho \geq 0}} 
\Bigg\{\int_0^T \int_\Omega \rho(x,t)\left(\frac{\partial \phi(x,t)}{\partial t} + \langle v(x,t), \nabla_x\phi(x,t)\rangle  - H^*\left(v(x,t),x, t\right) \right)dxdt \\- c\int_\Omega \phi(x,T)dx\Bigg\},
\end{split}
\end{equation*}
which is in the form of simple mean-field games. After applying a first-order monotone scheme with a numerical Hamiltonian $\hat H\colon (\Rn)^{2n} \times \Rn\times [0,T] \to \R$, the corresponding discretized saddle point formulation can be solved using PDHG. 
In~\cite{meng2023primal}, the 1-dimensional saddle point formula is discretized as 
\begin{equation}\label{eqt:saddle_pt_general_fulldisc}
\begin{aligned}
\min_{\substack{\phi_{i,k}\forall i,k\\ \phi_{i,1}=J(x_i)}}\max_{\substack{\rho_{i,k},v^+_{i,k}, v^-_{i,k}\forall i,k\\ \rho_{i,k} \geq 0}}  \sum_{i=1}^{n_x}\sum_{k=1}^{n_t-1} \rho_{i,k}\Bigg((D_t^- \phi)_{i,k+1}+ v^+_{i,k}(D_x^+\phi)_{i,k+1} +v^-_{i,k}(D_x^-\phi)_{i,k+1}  \\
- \hat H^*\left(v^+_{i,k}, v^-_{i,k}, x_i, t_{k+1}\right) \Bigg) - \frac{c}{\Delta t}\sum_{i=1}^{n_x}\phi_{i,n_t}.
\end{aligned}
\end{equation}
Here, \( x_1, \dots, x_{n_x} \) denote the discretization points of the spatial domain \([a,b]\), and \( t_1, \dots, t_{n_t} \) represent the discretization points of the temporal domain \([0,T]\). The values \( \phi_{i,k}, v^+_{i,k}, v^-_{i,k}, \rho_{i,k} \) correspond to the function values at the grid point \( (x_i, t_k) \).  
For finite difference approximations, we use \( D_t^- \), \( D_x^+ \), and \( D_x^- \) to denote the backward Euler temporal difference, right spatial difference, and left spatial difference, respectively. Specifically, these are defined as $(D_t^- \phi)_{i,k} = \frac{\phi_{i, k} - \phi_{i, k-1}}{\Delta t}$, $(D_x^+\phi)_{i,k}= \frac{\phi_{i+1,k} - \phi_{i,k}}{\Delta x}$, and $(D_x^-\phi)_{i,k} = \frac{\phi_{i,k} - \phi_{i-1,k}}{\Delta x}$.
According to the theory of first-order monotone schemes, the numerical Hamiltonian \( \hat{H} \) must satisfy two key properties:
\begin{itemize}
    \item Consistency: \( \hat{H}(p, p, x, t) = H(p, x, t) \).
    \item Monotonicity: \( \hat{H} \) must be non-increasing with respect to \( p^+ \) and non-decreasing with respect to \( p^- \).  
\end{itemize}
We denote the Fenchel-Legendre transform of \( \hat{H}(\cdot, \cdot, x, t) \) by \( \hat{H}^*(\cdot, \cdot, x, t) \). The PDHG method is then applied to compute the saddle point of problem~\eqref{eqt:saddle_pt_general_fulldisc}, where the first part of the solution provides the numerical approximation to the corresponding HJ PDE.
Fig.~\ref{fig:saddle_pt_results} presents a numerical result from~\cite{meng2023primal}, which solves the HJ PDE with the Hamiltonian \( H(p, x) = \|p\| + 3\exp(-4\|x - 1\|^2) + 1 \) and the initial condition \( J(x) = \sum_{i=1}^n \sin (\pi x_i) \).


In~\cite{meng2024primal}, the authors proposed solving the HJ PDE and the corresponding optimal control problem using the following saddle point formula
\begin{equation} \label{eqt:saddle_oc}
\begin{split}
\min_{\substack{\phi\\ \phi(x,0)=J(x)}}\max_{\rho\geq 0, \alpha} \Bigg\{ \int_0^T \int_\Omega \rho(x,t)\left(\frac{\partial \phi(x,t)}{\partial t} -\langle f_{x,t}(\alpha(x,t)), \nabla_x \phi(x,t)\rangle - L_{x,t}(\alpha(x,t)) \right)dxdt \\
- c\int_\Omega \phi(x,T)dx\Bigg\},
\end{split}
\end{equation}
where $f_{x,t}(\alpha)$ and $L_{x,t}(\alpha)$ denote $f(\alpha,x,T-t)$ and $L(\alpha,x,T-t)$, and $c$ is a positive hyperparameter. A one-dimensional discretization gives
\begin{equation*} 
\begin{split}
\min_{\substack{\phi\\ \phi_{i,1}=J(x_i)}}\max_{\rho\geq 0, \alpha}  \sum_{k=2}^{n_t} \sum_{i=1}^{n_x} \rho_{i,k}\Bigg((D_t^-\phi)_{i,k} - f_{i,k}(\alpha_{1,i,k})_+ (D_x^+\phi)_{i,k} - f_{i,k}(\alpha_{2,i,k})_- (D_x^-\phi)_{i,k} \\
- \hat L_{i,k}\left(\alpha_{1,i,k}, \alpha_{2,i,k}\right) \Bigg) - \frac{c}{\Delta t}\sum_{i=1}^{n_x} \phi_{i,n_t}.
\end{split}
\end{equation*}
Here, \( \hat{L}_{i,k} \) represents the numerical Lagrangian at \( (x_i, t_k) \), and other notations follow the conventions introduced earlier in~\cite{meng2023primal}. The resulting saddle-point problem is then solved using the PDHG method.  
Notably, this work can be seen as a generalization of~\cite{meng2023primal} for cases where \( f(\alpha, x, t) \neq \alpha \). When considering more general dynamics, the saddle-point formulation in~\eqref{eqt:saddle_oc} not only yields the solution \( \phi \) to the HJ PDE but also provides the optimal control function \( \alpha \) in a feedback form, making it particularly well-suited for optimal control problems.

\section{Conclusion}

This review explores several recent approaches for numerically solving Hamilton-Jacobi partial differential equations, ranging from classical grid-based methods to emerging deep learning techniques. 
Each approach offers distinct advantages: grid-based methods provide robust frameworks for handling complex geometries, representation formulas capture underlying mathematical structures, Laplace approximations provide sampling based schemes, deep learning methods show promise in handling complex boundary conditions and high-dimensional settings, and saddle point methods make connections to mean-field games and \update{introduce additional 
flexibility by allowing the use of primal–dual formulations, preconditioning strategies, 
and other optimization techniques.}
{Moreover, when high-dimensional HJ PDEs must be solved, one typically must resort to Lagrangian dynamics through optimal control formulations~\cite{ruthotto2020machine, lin2021alternating, bansal2021deepreach, kim2021hamilton, gelphman2025end}, often combined with deep learning techniques.}
While we highlight these recent developments, this review is by no means comprehensive. Notable approaches like optimal control algorithms~\cite{parkinson2020hamilton, onken2021neural, parkinson2018optimal, parkinson2022time}, optimal transport~\cite{onken2021ot, vidal2023taming, xu2023normalizing}, and mean-field games/control~\cite{ruthotto2020machine, chow2022numerical, carmona2021convergence, carmona2022convergence, vidal2024kernel, agrawal2022random, wang2025primal} offer alternative paths to solving HJ PDEs. 
{Some techniques designed for real-time control applications consider variations of the HJ PDE in~\eqref{eqt:HJ_general}, such as Hamilton-Jacobi-Isaacs equations~\cite{mitchell2005time, margellos2011hamilton, margellos2011hamilton, bansal2017hamilton, chen2021fastrack}.}
The field of numerical methods for Hamilton-Jacobi equations continues to evolve rapidly, with substantial work existing beyond what we have covered. 
{One key limitation of current approaches concerns the treatment of boundary conditions. Traditional numerical methods struggle in high dimensions, while neural network–based methods can only impose boundary conditions as soft constraints on sampled boundary points, except in cases where the conditions are trivial to enforce (e.g., box constraints). Nevertheless, new hybrid approaches continue to emerge, combining strengths from different methodological frameworks. The ongoing integration of these ideas, together with advances in computational capabilities, points to a promising future for the numerical solution of Hamilton–Jacobi PDEs.}

\section*{Acknowledge}
SWF was partially funded by NSF award DMS-2110745. S. Osher’s work was partially supported by ONR N00014-20-1-2787, NSF-2208272, STROBE NSF-1554564, and NSF 2345256. We would like to thank Dr. Mo Zhou for the discussion on Section 5.1.

\bibliographystyle{abbrv}
\bibliography{references}

\end{document}